\documentclass[12pt,reqno]{amsart}
\usepackage{mathrsfs}
\usepackage{amsgen}
\usepackage{amscd}
\usepackage{amsmath}
\usepackage{latexsym}
\usepackage{amsfonts}
\usepackage{amssymb}
\usepackage{amsthm}
\usepackage{graphicx}
\usepackage{color}
\usepackage{graphicx}
\usepackage{enumerate}
\usepackage{amssymb,amsmath,graphicx,amsfonts,euscript}
\usepackage{mathrsfs}
\usepackage[hyphens]{url}
\usepackage{caption}
\usepackage{subcaption}
\usepackage[margin=3cm, a4paper]{geometry}
\DeclareGraphicsExtensions{.eps}

\newtheorem{theorem}{Theorem}[section]
\newtheorem{lemma}[theorem]{Lemma}

\theoremstyle{definition}

\newtheorem{example}[theorem]{Example}
\newtheorem{assumption}[theorem]{Assumption}

\theoremstyle{remark}

\def\R{\mathbb{R}}

\def\N{\mathbb{N}}
\def\L{\mathcal{L}}

\newcommand{\s}{\mathop{\mathrm{span}}}
\newcommand{\re}{\mathop{\mathrm{Re}}}
\newcommand{\Arg}{\mathop{\mathrm{arg}}}

\newcommand{\black}{\color{black}}

\numberwithin{equation}{section}

\allowdisplaybreaks

\begin{document}

\title[]{EXPLICIT EXPONENTIAL RUNGE-KUTTA METHODS FOR SEMILINEAR INTEGRO-DIFFERENTIAL EQUATIONS}

\author[A. Ostermann]{Alexander Ostermann}
\address{\hspace*{-12pt}A.~Ostermann: Department of Mathematics, Universit\"{a}t Innsbruck, Technikerstrasse 13, 6020 Innsbruck, Austria}
\email{alexander.ostermann@uibk.ac.at}

\author[F. Saedpanah]{Fardin Saedpanah}
\address{\hspace*{-12pt}F.~Saedpanah: Department of Mathematics, University of Kurdistan, P. O. Box 416, Sanandaj, Iran and Department of Engineering, University of Bor{\aa}s, SE-501 90 Bor{\aa}s, Sweden}
\email{f.saedpanah@uok.ac.ir; fardin.saedpanah@hb.se}

\author[N. Vaisi]{Nasrin Vaisi}
\address{\hspace*{-12pt}N.~Vaisi: Department of Mathematics, University of Kurdistan, P. O. Box 416, Sanandaj, Iran}
\email{nasrin\_vaisi@yahoo.com; n.vaisy@sci.uok.ac.ir}



\begin{abstract}\noindent
The aim of this paper is to construct and analyze explicit exponential Runge--Kutta
 methods for the temporal discretization of linear and semilinear integro-differential
 equations. By expanding the errors of the numerical method in terms of the solution, we derive
 order conditions that form the basis of our error bounds for integro-differential equations.
 The order conditions are further used for constructing numerical methods. The convergence
 analysis is performed in a Hilbert space setting, where the smoothing effect of the resolvent
 family is heavily used. For the linear case, we derive the order  conditions for
 general order $p$ and prove  convergence of order $p$, whenever these conditions are satisfied.
 In the semilinear case, we consider in addition spatial discretization by a spectral Galerkin
 method, and we  require locally Lipschitz continuous nonlinearities. We derive
 the order conditions for orders one and two, construct methods satisfying these conditions
 and prove their  convergence. Finally, some numerical experiments illustrating our theoretical
 results are given.
\end{abstract}

\maketitle

\section{Introduction}\label{sec:introduction}
In this paper we consider the time discretization of linear integro-differential equations
\begin{eqnarray}\label{linear integro-differential equation}
\frac{\partial u(x,t)}{\partial t} +\int_0^t b(t-s)Au(x,s) ds = f(x,t), \qquad u(x,0)=u_0(x),
\end{eqnarray}
 and the full discretization of semilinear integro-differential equations of the form
\begin{eqnarray} \label{semilinear integro-differential equation}
 \frac{\partial u(x,t)}{\partial t} +\int_{0}^{t}b(t-s)Au(x,s)\,ds=f(x,t,u(x,t)), \qquad  u(x,0)=u_{0}(x),
\end{eqnarray}
for $x$ in a domain $\Omega\subseteq \R^{d}$ and $t \in [0,T]$, taken together with homogeneous Dirichlet boundary
 conditions. The operator $A$ is self-adjoint and positive definite on a Hilbert space $H$
 with compact inverse. The kernel $b$ is assumed  to be real-valued and positive definite,
 i.e., for each $T >0$ the kernel $b$ belongs to $L^{1}(0,T)$ and satisfies
 \[ \int_{0}^{T} \psi(t) \int_{0}^{t}b(t-s)\psi(s)\,ds\,dt\geqslant 0 \quad {\text{for all}} \quad \psi \in C[0,T].\]

\black Semilinear problems, or linear versions thereof, \black are used to model viscoelasticity and heat
 conduction in materials with memory, see, e.g., \cite{CLP, KLS, LST, MP, MT}.
\black When the kernel $b$ is weakly singular, one can interpret the evolution equation as a fractional wave equation, see \cite{ Garrappa}. When \black the kernel $b$ is smooth such equations are hyperbolic in nature, \black while when \black $b$ has a weak
 singularity at $t = 0$, they exhibit certain features of parabolic equations. As a typical weakly singular example,
 we mention the Riesz kernel
\[b(t)=\dfrac{t^{\beta-1}}{\Gamma(\beta)}, \quad 0<\beta<1.\]

We recall that $b$ is positive definite if and only if
$$ \re \black\mathcal{L}(b)(i\theta)\black=\int_{0}^{\infty}b(t)\cos(\theta t)\,dt \geqslant 0 \quad {\text{for all}} \quad \theta \in \mathbb{R},$$
 where $\mathcal{L}(b)$ denotes the Laplace transform of $b$. A sufficient condition for this
 to hold is that $b \in L_\text{loc}^{1} \cap C^{2}(0,\infty)$, $(-1)^{n}b^{(n)}\geqslant 0$ for all $t>0$, $n=0, 1, 2$, and that
$b^{(2)}(t)$ is nonincreasing and convex, i.e., $b$ is 4-monotone kernel, see \cite[Definition 3.4]{Pruss}.

 The numerical solution of problem~(\ref{semilinear integro-differential equation})  has been studied in, e.g., \cite{CLP, CP, CT, KLS, LST, MT, Sanz}.
 The methods considered in \cite{LST, MT} are based on the finite element method for the spatial discretization,
 together with the first- and second-order backward difference methods or the Crank--Nicolson method in time, with appropriate quadrature  formulas applied to the convolution term.
In \cite{CLP}, by considering the Riesz kernel, a systematic and computationally affordable approach was derived.
 It gives second order accuracy in time under realistic regularity assumptions.

For differential equations, the idea of exponential integrators is an old one and has been proposed independently  by many authors. The numerical comparisons presented in \cite{KT, K} show a number of examples for which explicit exponential integrators perform better than standard integrators. \black In particular, exponential integrators provide exact solutions for linear homogeneous problems, and high-order approximations to linear inhomogeneous problems. As a consequence, very accurate numerical solutions can be obtained with large time steps even for nonsmooth and weakly singular kernels, which is an issue in integro-differential equations. The convergence behavior of implicit and linearly implicit Runge--Kutta methods  for parabolic problems was studied in \cite{LO1, LO2}, that of implicit exponential Runge--Kutta methods in \cite{HO3}. Later, in a series of papers, \black new techniques were introduced for proving error bounds in the explicit case. In \cite{HO1, HO2} the authors derived the order conditions for
 stiff problems and, based on these, proved error bounds for parabolic problems.
 The new conditions enabled them to analyze the methods presented in the literature and, in addition,
 to develop new methods that do not suffer from reduced orders. In \cite{KLS}
 the exponential Euler method was generalized to a stochastic version of these problems.
  The resulting scheme was named Mittag-Leffler--Euler integrator. Our aim with this paper
 is to give error bounds for the time discretization of integro-differential equations by exponential Runge--Kutta
 methods. A fully discrete scheme is then obtained by combining the time discretization with the spectral Galerkin method for spatial discretization.

The outline of the paper is as follows. After presenting the abstract framework and some preliminaries,
 we state our main assumptions for the linear problem in Section~\ref{Preliminaries}. Then, in
Section \ref{Linear_problems}, we study linear problems and introduce our numerical scheme for the temporal semidiscretization,
 viz.~(\ref{eq:EQR}). In Theorem~\ref{th.1}, we state and prove the convergence result for exponential
 Runge--Kutta methods. In Section \ref{Semilinear_problems}, we define a general
 class of exponential Runge--Kutta methods for semilinear integro-differential equations, and introduce the fully discrete scheme.
 Our main results are contained in Section \ref{Convergence}, where we derive order conditions for explicit exponential
 Runge--Kutta methods of order two applied to semilinear problems.
 For the analysis of~(\ref{semilinear integro-differential equation}), an abstract Hilbert space framework of
 locally Lipschitz continuous nonlinearities is chosen and the smoothing effect of the resolvent is used.
 Based on the order conditions, we obtain explicit exponential Runge--Kutta methods of order
 two and show their convergence. The convergence results for the exponential Euler method and for second-order
 methods are given in Theorems \ref{th.2} and \ref{th.3} respectively. Finally, in Section \ref{Numerical}, we
 present some numerical experiments which illustrate our theoretical results.
\section{Preliminaries}\label{Preliminaries}
\subsection{The abstract setting.}
Let $H$ be a real, separable, infinite-dimensional Hilbert space.  An important example is $H=L^{2}(\mathcal{D})$.
The standard inner product and norm in $H$ will be denoted by $\|\cdot\|$ and $(\cdot ,\cdot)$, respectively.
The space of all bounded linear operators on $H$ will be denoted by \black $\mathcal{B}=\mathcal{B}(H)$\black.

\begin{assumption}\label{A_operator}
 Let $A$ be a self-adjoint, positive definite operator on the Hilbert space $H$ with compact inverse, and let the kernel $b$
 be positive definite.
\end{assumption}

\noindent
The standard example is $A =-\Delta$  with homogeneous
 Dirichlet boundary conditions on an open and bounded domain $\mathcal{D}\subseteq \R^{d}$. This operator is positive definite
 on $L^{2}(\mathcal{D})$ with an orthonormal eigenbasis $\{\psi_{j}\}_{j=1}^{\infty}$ and corresponding
 eigenvalues $\{\lambda_{j}\}_{j=1}^{\infty}$ such that
\begin{eqnarray*}
 A\psi_{j}=\lambda_{j}\psi_{j},\quad 0<\lambda_{1}\leq \lambda_{2}\leq \cdots \leq \lambda_{j}\leq \cdots ,\quad  \lambda_{j}\rightarrow \infty.
\end{eqnarray*}

\subsection{Resolvent family} Under Assumption~\ref{A_operator} it follows from \cite[Corollary 1.2]{Pruss} that there
 exists a strongly continuous family $\{S(t)\}_{t\geq 0}$ of bounded linear operators on $H$ such
 that the function  $u(t)=S(t)u_{0}$, $u_{0}\in H$, is the unique solution of
\[u(t)+A\int_{0}^{t}B(t-s)u(s)\,ds=u_{0},\quad t\geqslant 0,\]
with $B(t)=\int_{0}^{t}b(s)\,ds.$ If $t\mapsto u(t)=S(t)u_{0}$ is differentiable for $t > 0$, then $u$ is
 the unique solution of
\[ u'(t)+A\int_{0}^{t}b(t-s)u(s)\,ds=0,\quad t>0,\quad u(0)=u_{0}.\]
We refer to the monograph \cite{Pruss} for a comprehensive theory of resolvent families for
 Volterra equations.  An important feature of the resolvent family $\{S(t)\}_{t\geq 0}$ is that
 it does not have the semigroup property; that is, $S(t + s)\neq S(t)S(s)$, in general. This is the
 mathematical reflection of the fact that the solution possesses a nontrivial memory. In
 our special setting, using the spectral decomposition of $A$, an explicit representation
 of $S(t)$ is given by the Fourier series
\begin{eqnarray}\label{explicit_representation}
S(t)v=\sum_{k=1}^{\infty}s_{k}(t)(v,\psi_{k})\psi_{k},
\end{eqnarray}
where the functions $s_{k}(t)$ are the solutions of the \black ordinary integro-differential equations\black
\begin{eqnarray}\label{equation}
s'_{k}(t)+\lambda_{k}\int_{0}^{t}b(t-s)s_{k}(s)\,ds=0,\quad t>0,\quad s_{k}(0)=1,
\end{eqnarray}
with $\lbrace(\lambda_{k},\psi_{k})\rbrace_{k=1}^{\infty}$ being the eigenpairs of $A$.

The following assumption, which establishes the smoothing property of the
 resolvent family $\{S(t)\}_{t\geq 0}$, is one of the central tools for proving the main results of this paper.

 \begin{assumption}\label{smooth_property}
  We assume that the resolvent family $\{S(t)\}_{t\geq 0}$ is strongly continuous for $t\ge 0$ and strongly continuously
 differentiable for $t > 0$ and enjoys the following smoothing property:
 there are constants $C$ and $1<\rho <2$ such that for any $0<t\le T$, we have
\begin{equation}\label{property}
\|A^{\alpha}S(t)\|_{\mathcal{B}}\leq Ct^{-\alpha\rho},\qquad 0\le \alpha < \tfrac1{\rho} .
\end{equation}
\end{assumption}

\noindent
 This smoothing property is verified in \cite[Theorem 5.5]{MT} for the Riesz kernel
 ${t^{\beta-1}}/{\Gamma(\beta)}$, $0<\beta<1$,  with $\rho=\beta+1$. A more general class of
 kernels $b$ for which~(\ref{property}) is satisfied is  the class of 4-monotone kernels with
\[ \rho=1+\frac{2}{\pi}\sup\lbrace\lvert\Arg \L(b)(z)\rvert, \re z>0 \rbrace \in (1,2), \]
and $\L(b)(z)\leq Cz^{1-\rho}$ for $z>1,$ where this latter condition may be substituted by the
condition \black$|b(t)|\leq Ct^{\rho-1}$, $t\in (0,1)$\black, see \cite[Remarks 2.5, 3.8 and Lemma A.4]{BGK}. In
particular, $b$ does not have to be analytic.
\section{linear problems: exponential quadrature}\label{Linear_problems}
In this section, we derive error bounds for exponential Runge--Kutta discretizations
of linear integro-differential equations~(\ref{linear integro-differential equation})
with a time-invariant operator $A$, $u_{0}\in H.$
 We consider problems with $f$ being smooth, so that we can expand the solution in a Taylor series.

 Let $ 0=t_{0} <t_{1} < \cdots <t_{M}=T$ be a uniform partition of the time interval $\left[ 0,T\right]$
 with time step $ h=t_{m+1}-t_{m}$, $m=0,1,\dots,M-1.$ Under Assumption~\ref{A_operator} there exists a resolvent family $\{S(t)\}_{t\geq 0}$
 of bounded linear operators on $H$, which is strongly continuous for $t\geq 0$ and differentiable for $t>0$, such that for $m=0,1,\dots,M,$ by using the variation-of-constants formula, we have
    \begin{align}\label{exact_sol1}
    u(t_{m})&=S(t_{m})u_{0}+\int^{t_{m}}_{0}S(t_{m}-\sigma)f(\sigma)\,d\sigma\nonumber \\
                  &=S(t_{m})u_{0}+\sum_{j=0}^{m-1}\int^{h}_{0}S(t_{m-j}-\sigma)f(t_{j}+\sigma)\,d\sigma.
    \end{align}
 A scheme is obtained by approximating the function $f$ within the integral by its
 interpolation polynomial, using the quadrature nodes \black $0=c_{1}<c_{2}<\cdots< c_{s}\leq 1$\black. This yields an \textit{exponential quadrature rule}, for $m=0,1,\dots,M,$
\begin{subequations}\label{numerical_sol1}
\begin{eqnarray}\label{a}
 U_{m}=S(t_{m})u_{0}+h\sum_{j=0}^{m-1}\sum_{i=1}^{s}b_{i}(t_{m-j})f(t_{j}+c_{i}h),
  \end{eqnarray}
  with weights
 \begin{eqnarray}\label{b}
  b_{i}(t_{l})=\frac{1}{h}\int^{h}_{0}S(t_{l}-\sigma)L_{i}(\sigma)\, d\sigma,\quad 1\leqslant l \leqslant m,
   \end{eqnarray}
   \end{subequations}
  where $L_{i}$ are the Lagrange interpolation polynomials
  \[
  L_{i}(\sigma)=\prod_{n=1,\;n\neq i}^{s}\dfrac{\sigma/h -c_{n}}{c_{i}-c_{n}},\quad i=1,\ldots,s.
  \]
  We need the weights $b_{i}(t_l)$ to be uniformly bounded in $h \geqslant 0$.
  Since the weights $b_{i}(t_l)$ are linear combination of the operators
  \begin{eqnarray}\label{phi}
  \varphi_{k,h}(t_{l})=\frac{1}{h^{k}}\int^{h}_{0}S(t_{l}-\sigma)\dfrac{\sigma^{k-1}}{(k-1)!}\, d\sigma,\quad k\geqslant 1,
  \end{eqnarray}
   we will use the following important lemma.
  \begin{lemma}  Under Assumption \ref{smooth_property}, the operators $\varphi_{k,h}(t_{l})$,
  $1\leqslant l\leqslant M$, $k\geqslant 1,$  are bounded on $H$.
  \end{lemma}
  \begin{proof}
  The estimate of $\varphi_{k,h}(t_{l})$ is a consequence of~(\ref{property}) with $\alpha =0$, as
  \begin{align*}\black
 \| \varphi_{k,h}(t_{l}) \|_{\mathcal{B}}\leqslant \frac{1}{h^{k}}\int^{h}_{0} \|S(t_{l}-\sigma)\|_{\mathcal{B}}\dfrac{\sigma^{k-1}}{(k-1)!}\, d\sigma\leqslant \frac{C}{k!} \black
  \end{align*}
   is obviously bounded (uniformly for $h>0$).
  \end{proof}
   Therefore, for the coefficients of the exponential Runge--Kutta method~(\ref{numerical_sol1}), we get a
 smoothing property similar to~(\ref{property}), that is, for given $1<\rho <2$,
 \begin{eqnarray}\label{property1}
\black \|A^{\alpha}\phi(t_l)\|_{\mathcal{B}}\leq Ct_l^{-\alpha\rho},\black\quad 0\le \alpha < 1/\rho,
\end{eqnarray}
for $\phi =b_{i}$, $i=1,\ldots, s$, \black and $1\le l\le M$.\black

\black Exponential quadrature rules for linear integro-differential equations can also be formulated from scratch in the following way. For $1\le i\le s$ and $1\le l\le M$, let $b_i(t_l)$ denote bounded operators (with a bound that is uniform in the step size $h$). For nonconfluent nodes $0=c_1<c_2 <\ldots<c_s$, we consider the following exponential quadrature rule for the time discretization of \eqref{exact_sol1}:
\begin{equation}\label{eq:EQR}
 U_{m}=S(t_{m})u_{0}+h\sum_{j=0}^{m-1}\sum_{i=1}^{s}b_{i}(t_{m-j})f(t_{j}+c_{i}h),\qquad 0\le m\le M.
\end{equation}
The weights $b_i(t_l)$ and the nodes $c_i$ have to satisfy certain order conditions, which will be studied next.\black
  \subsection{Error expansion and order conditions.}
In order to analyze~(\ref{eq:EQR}), we expand the exact solution~(\ref{exact_sol1}) into a Taylor series with
 remainder in integral form
\begin{align}\label{Taylor series exact}
  u(t_{m})&=S(t_{m})u_{0}+\sum_{j=0}^{m-1}\int^{h}_{0}S(t_{m-j}-\sigma)f(t_{j}+\sigma)\,d\sigma\nonumber\\
                &=S(t_{m})u_{0}+\sum_{j=0}^{m-1}\int^{h}_{0}S(t_{m-j}-\sigma)\sum_{k=0}^{p-1}\frac{\sigma^k}{k!}f^{(k)}(t_{j})\,d\sigma\\
                &\quad +\sum_{j=0}^{m-1}\int^{h}_{0}S(t_{m-j}-\sigma)\int_{0}^{\sigma}\frac{(\sigma-\tau)^{p-1}}{(p-1)!}f^{(p)}(t_{j}+\tau)\,d\tau \,d\sigma.\nonumber
\end{align}
Now this is compared with the Taylor series of the numerical solution~(\ref{eq:EQR})
\begin{align}\label{Taylor series numerical}
  U_{m}&=S(t_{m})u_{0}+h\sum_{j=0}^{m-1}\sum_{i=1}^{s}b_{i}(t_{m-j})f(t_{j}+c_{i}h)\nonumber\\
                &=S(t_{m})u_{0}+h\sum_{j=0}^{m-1}\sum_{i=1}^{s}b_{i}(t_{m-j})\sum_{k=0}^{p-1}\frac{c_{i}^{k}h^{k}}{k!}f^{(k)}(t_{j})\\
                &\quad +h\sum_{j=0}^{m-1}\sum_{i=1}^{s}b_{i}(t_{m-j})\int_{0}^{c_{i}h}\frac{(c_{i}h-\sigma)^{p-1}}{(p-1)!}f^{(p)}(t_{j}+\sigma)\,d\sigma.\nonumber
\end{align}
 By subtracting~(\ref{Taylor series numerical}) from~(\ref{Taylor series exact}), we get
 \begin{align}\label{linear error}
 u(t_{m})-U_{m}&=\sum_{j=0}^{m-1}\sum_{k=0}^{p-1}\frac{1}{k!}\Big(\int^{h}_{0}S(t_{m-j}-\sigma)\sigma^{k} \,d\sigma-h^{k+1}\sum_{i=1}^{s}c_{i}^{k}b_{i}(t_{m-j})\Big)f^{(k)}(t_{j})\nonumber\\
                          &\quad +\sum_{j=0}^{m-1}\int^{h}_{0}S(t_{m-j}-\sigma)\int_{0}^{\sigma}\frac{(\sigma-\tau)^{p-1}}        {(p-1)!}f^{(p)}(t_{j}+\tau)\,d\tau \,d\sigma\\
                         &\quad -h\sum_{j=0}^{m-1}\sum_{i=1}^{s}b_{i}(t_{m-j})\int_{0}^{c_{i}h}\frac{(c_{i}h-\sigma)^{p-1}}{(p-1)!}f^{(p)}  (t_{j}+\sigma)\,d\sigma. \nonumber
 \end{align}
 The coefficients
 \begin{equation}\label{order-condition}
 M_{k}(t_{l})=\int^{h}_{0}S(t_{l}-\sigma)\sigma^{k-1}\,d\sigma-h^{k}\sum_{i=1}^{s}c_{i}^{k-1}b_{i}(t_{l}),\quad 1\leqslant l\leqslant m,\quad k=1,\ldots,p,
 \end{equation}
  of the low-order terms in~(\ref{linear error}) being zero turn out to be the \textit{order conditions }of the exponential Runge--Kutta method~(\ref{eq:EQR}).
  The order conditions for an $s$-stage exponential quadrature rule are given in Table~\ref{tab1}. \black It is easy to verify that the method~\eqref{numerical_sol1} satisfies these conditions up to order $p=s$.\black

 \begin{table}[h]
 \caption{Order conditions for an $s$-stage exponential quadrature rule~\eqref{eq:EQR}. The functions $\varphi_{k,h}$ are defined in~\eqref{phi}.\label{tab1}}
 \renewcommand{\arraystretch}{1.8}
 \begin{tabular}{|c|c|}
 \hline
 $\text{Order}$ & $\text{Order condition}$ \\
 \hline
 $1$ & $ \sum_{i=1}^{s}b_{i}(t_{l})=\varphi_{1,h}(t_{l})$ \\
\hline
 $2$ & $ \sum_{i=1}^{s}b_{i}(t_{l})c_{i}=\varphi_{2,h}(t_{l})$ \\
\hline
 $\vdots $ & $ \vdots $ \\
\hline
 $ p $ & $\rule[-14pt]{0pt}{36pt}\sum_{i=1}^{s}b_{i}(t_{l})\dfrac{c_{i}^{p-1}}{(p-1)!}=\varphi_{p,h}(t_{l})$\\
\hline
\end{tabular}
\end{table}

 An exponential quadrature method has order $p$, if all conditions in Table~\ref{tab1} are satisfied. Note that the order conditions are linear in the weight functions
$b_i(t_l)$ and form a Vandermonde system for given pairwise distinct nodes $c_1,\ldots,c_s$. Therefore, by choosing $s=p$, the weights $b_i(t_l)$ of an $s$-stage exponential quadrature rule of order $p=s$ are uniquely defined in terms of the given nodes.

We are now ready to state our convergence result.
\begin{theorem}\label{th.1}
Let Assumptions~\ref{A_operator} and~\ref{smooth_property} be satisfied. For the numerical solution of~(\ref{linear integro-differential equation}), we consider an exponential Runge--Kutta method~(\ref{eq:EQR}) of order $p\ge 1$. If $f^{(p)}\in L^{1}(0,T)$ then the following error bound holds
\begin{equation*}
\|u(t_{m})-U_{m}\|\leqslant C h^{p} \int_{t_{0}}^{t_{m}}\| f^{(p)}(\tau)\| \,d\tau,
\end{equation*}
uniformly on $0\leqslant t_{m} \leqslant T$. The constant $C$ depends on the final time $T$, but is independent of $m$ and $h$.
\end{theorem}
\begin{proof}
Using the order conditions $M_{k}(t_{l})=0$ in~(\ref{linear error}), we have
 \begin{align*}
 u(t_{m})-U_{m}&=\sum_{j=0}^{m-1}\int^{h}_{0}S(t_{m-j}-\sigma)\int_{0}^{\sigma}\frac{(\sigma-\tau)^{p-1}}{(p-1)!}f^{(p)}(t_{j}+\tau) \,d\tau \,d\sigma\\
                         &\quad -h\sum_{j=0}^{m-1}\sum_{i=1}^{s}b_{i}(t_{m-j})\int_{0}^{c_{i}h}\frac{(c_{i}h-\sigma)^{p-1}}{(p-1)!}f^{(p)}  (t_{j}+\sigma)\,d\sigma. \nonumber
 \end{align*}
By changing the order of integration, taking norms, and using the smoothing properties~(\ref{property}) and~(\ref{property1}) with $\alpha=0$, we obtain
\begin{align*}
\| u(t_{m})-U_{m}\| & \leqslant C\sum_{j=0}^{m-1}\int^{h}_{0}\int_{0}^{\sigma} \frac{(\sigma-\tau)^{p-1}}{(p-1)!} \|f^{(p)}(t_{j}+\tau)\| \,d\tau \,d\sigma\\
                                & \qquad \black+Ch\sum_{j=0}^{m-1}\max_{1\leqslant i\leqslant s}\int_{0}^{c_{i}h}\frac{h^{p-1}}{(p-1)!}\|f^{(p)}(t_{j}+\sigma)\| \,d\sigma\\ \black
                                & \leqslant C\sum_{j=0}^{m-1}\int^{h}_{0}\|f^{(p)}(t_{j}+\tau)\| \int_{\tau}^{h}
                                \frac{(\sigma-\tau)^{p-1}}{(p-1)!} \,d\sigma \,d\tau\\
                                & \qquad +\black Ch\sum_{j=0}^{m-1}\int_{t_{j}}^{t_{j+1}}\frac{h^{p-1}}{(p-1)!}\|f^{(p)}(\sigma)\| \,d\sigma\\ \black
                               & \leqslant Ch^{p} \int_{t_{0}}^{t_{m}}\| f^{(p)}(\tau)\| \,d\tau.
\end{align*}
This is the desired result.
\end{proof}
  \section{Semilinear problems: exponential Runge--Kutta methods}\label{Semilinear_problems}
For the numerical solution of semilinear problems~(\ref{semilinear integro-differential equation}), we proceed analogously
to the construction of exponential Runge--Kutta methods for differential equations. We start from the variation-of-constants formula
\begin{align}\label{V_O_C}
 u(t_{m})&=S(t_{m})u_{0}+\int^{t_{m}}_{0}S(t_{m}-\sigma)f(\sigma, u(\sigma))\,d\sigma \nonumber\\
               &=S(t_{m})u_{0}+\sum_{j=0}^{m-1}\int^{h}_{0}S(t_{m-j}-\sigma)f(t_{j}+\sigma, u(t_{j}+\sigma))\,d\sigma.
\end{align}
Here $\{S(t)\}_{t\geq 0}$ is a resolvent family of bounded linear operators on $H$, which is
 strongly continuous for $t\geq 0$ and differentiable for $t>0$, $u_{0}\in H $, $f\in L^{\infty}([0,T];H)$.
  We note that the resolvent family does not enjoy the semigroup property due to the
nonlocality of the  kernel in~(\ref{semilinear integro-differential equation}).

 The numerical scheme is defined recursively for $m\ge 1$ by
\begin{subequations}\label{(1.9)}
\begin{align}
&U_{m}=S(t_{m})u_{0}+h\sum^{m-1}_{j=0}\sum^{s}_{i=1}b_{i}(t_{m-j})f(t_{j}+c_{i}h,U_{j,i})\label{(1.9a)}\\
\noalign{\noindent and}
& U_{m-1,q}=S(t_{m-1}+c_{q}h)u_{0} + h\sum_{k=1}^{q-1}a_{qk}f(t_{m-1}+c_{k}h,U_{m-1,k})\nonumber\\
                  & \qquad \qquad+h\sum^{m-2}_{l=0}\sum^{s}_{i=1}b_{i}^{q}(t_{m-l-1})f(t_{l}+c_{i}h,U_{l,i}), \quad 1\leqslant q \leqslant s, \label{(1.9b)}
\end{align}
\end{subequations}
where $U_{m}$ denotes the numerical approximation to $u(t_{m})$ and $U_{m-1,q}\approx u(t_{m-1}+c_qh)$. Here,
 the method's coefficients $a_{qk}$, $b_{i}^{q}$ and $b_{i}$ are constructed from the resolvent family $\{S(t)\}_{t\geq 0}$, in general. Therefore, it is plain to assume that the coefficients satisfy a smoothing property similar to~(\ref{property1}) for $\phi=b_{i}$, $\phi=b_{i}^{q}$ and $\phi=a_{qk}$, for $i,q=1,\ldots, s$ and $k=1,\ldots, q-1$.

The above scheme is called an \textit{explicit exponential Runge--Kutta method} for integro-differential equations. \black As in the linear case, we will always assume that $c_1=0$.\black

\subsection{Discretisation in space.}
For spatial discretization, we define a finite dimensional subspace $H_{N}$ of $H$ by
 $H_{N}= \s \lbrace\psi_{1},\cdots,\psi_{N}\rbrace$, where $\lbrace\psi_{k}\rbrace_{k=1}^{\infty}$ are the eigenvectors of $A$, i.e., $A\psi_{k}=\lambda_{k}\psi_{k},~ k\in \mathbb{N}$.
  Further, we define the projector
\begin{eqnarray}\label{Projector}
\mathcal{P}_{N}:H\rightarrow H_{N},\quad \mathcal{P}_{N}v=\sum_{k=1}^{N}(v,\psi_{k})\psi_{k},\quad v \in H.
\end{eqnarray}
We also consider the projected operator
\begin{eqnarray}\label{Operator}
A_{N}:H_{N}\rightarrow H_{N},~~A_{N}=A\mathcal{P}_{N},
\end{eqnarray}
which generates a family of resolvent operators $ \lbrace S_{N}(t)\rbrace_{t\geqslant0}$ in $H_{N}$. It is clear that
\begin{eqnarray}\label{Resolvent_operators}
S_{N}(t)\mathcal{P}_{N}=S(t)\mathcal{P}_{N},
\end{eqnarray}
 and also
\begin{align*}
\Vert A^{-\nu}(I-\mathcal{P}_{N})x\Vert^{2}&=\sum_{k=1}^{\infty}\lambda_{k}^{-2\nu}((I-\mathcal{P}_{N})x,\psi_{k})^{2}
                                                                     =\sum_{k=N+1}^{\infty}\lambda_{k}^{-2\nu}(x,\psi_{k})^{2},\\
                                                                     &\leq \sup_{k\geqslant N+1} \lambda_{k}^{-2\nu}\sum_{k=N+1}^{\infty}(\psi_{k})^{2}
                                                                     \leq \lambda_{N+1}^{-2\nu}\sum_{k=N+1}^{\infty}(x,\psi_{k})^{2}
                                                                     \leq \lambda_{N+1}^{-2\nu}\|x\|^{2}.
\end{align*}
So
\begin{eqnarray}\label{Bound}
\Vert A^{-\nu}(I-\mathcal{P}_{N})\Vert=\sup_{k\geqslant N+1}\lambda_{k}^{-\nu}=\lambda_{N+1}^{-\nu},\quad \nu\geqslant 0.
\end{eqnarray}
The representation of $S_{N}$, similar to~(\ref{explicit_representation}), is given by
\[ S_{N}(t)v=\sum_{k=1}^{N}s_{k}(t)(v,\psi_{k})\psi_{k}. \]
This motivates us to consider the following fully discrete approximation of~(\ref{semilinear integro-differential equation}), based on the temporal
approximation~(\ref{(1.9)}):
\begin{subequations}\label{1.14}
\begin{align}
& U_{m}^{N}=S_{N}(t_{m})\mathcal{P}_{N}u_{0}+h\sum^{m-1}_{j=0}\sum^{s}_{i=1}b_{i}^{N}(t_{m-j}) \mathcal{P}_{N} f(t_{j}+c_{i}h,U_{j,i}^{N}),\label{1.14a}\\
& U_{m-1,q}^{N}=S_{N}(t_{m-1}+c_{q}h)\mathcal{P}_{N}u_{0} + h\sum_{k=1}^{q-1}a_{qk}^{N}\mathcal{P}_{N} f(t_{m-1}+c_{k}h,U_{m-1,k}^{N})\nonumber\\
                  & \qquad \qquad+h\sum^{m-2}_{l=0}\sum^{s}_{i=1}b_{i}^{q,N}(t_{m-l-1})\mathcal{P}_{N} f(t_{l}+c_{i}h ,U_{l,i}^{N}),\qquad 1\le q\le s,\label{1.14b}
\end{align}
\end{subequations}
 where the coefficients $b_{i}^N(t)$, $b_{i}^{q,N}(t)$ and $a_{qk}^N$ are simply given by
\begin{eqnarray*}
b_{i}^{N}(t)=\mathcal{P}_{N}b_{i}(t),\qquad b_{i}^{q,N}(t)=\mathcal{P}_{N}b_{i}^{q}(t),\qquad a_{qk}^{N}=\mathcal{P}_{N}a_{qk}.
\end{eqnarray*}
They are bounded operators on $H_N$ and satisfy a smoothing property similar to~(\ref{property1}), but now uniformly in $N\in\N$. In this paper, due to the particular choice of the coefficients in \eqref{(1.9)}, the relations
\begin{eqnarray}\label{coefficient}
b_{i}^{N}(t)\mathcal{P}_{N}=b_{i}(t)\mathcal{P}_{N},\qquad b_{i}^{q,N}(t)\mathcal{P}_{N}=b_{i}^{q}(t)\mathcal{P}_{N},\qquad a_{qk}^{N}\mathcal{P}_{N}=a_{qk}\mathcal{P}_{N}
\end{eqnarray}
will always hold. By spectral theory we also define $V=\mathcal{D}(A^{\nu})$ with norm
$$\|v\|_{V}^{2}=\|A^{\nu}v\|^{2}=\sum_{k=1}^{\infty}\lambda_{k}^{\nu}(v,\psi_{k})^{2},\quad  \nu \in \mathbb{R},~ v\in V.$$
Our main assumptions on the nonlinearity $f$ are those of \cite{Henry, Pazy}. In particular, we make the following assumption.
\begin{assumption}\label{assumption1}
Let $f: [0, T] \times V \rightarrow H$ be locally Lipschitz continuous in a strip along
 the exact solution $u$. Thus there exists a real number $L(R, T)$ such that
\[\|f(t, v) - f(t,w)\| \leqslant L\|v - w\|_{V}\]
for all $t \in [0, T]$ and $\max(\|v - u(t)\|_{V}, \|w - u(t)\|_{V}) \leqslant R$.
\end{assumption}

\section{Convergence results for semilinear problems}\label{Convergence}
We are now in a position to prove the convergence properties of exponential Runge--Kutta
methods for the semilinear problem~(\ref{semilinear integro-differential equation}). For simplicity in presentation, we limit our analysis to
methods of orders one and two.
\subsection{Convergence of the exponential Euler integrator.}For $s = 1$, the only
 reasonable selection is the exponential form of Euler's method with $b_1(t_l) = \varphi_{1,h}(t_l)$ \black and $c_1=0$. It will be called exponential Euler integrator. \black
 Applied to the space discretization of~(\ref{semilinear integro-differential equation}), it has the form
\begin{eqnarray}\label{numerical}
 U^{N}_{m}=S_{N}(t_{m})\mathcal{P}_{N} u_{0}+h\sum^{m-1}_{j=0}b_{1}^{N}(t_{m-j})\mathcal{P}_{N}f(t_{j},U_{j}^{N}),
\end{eqnarray}
with
\[b_{1}^{N}(t_{l})=\frac{1}{h}\int_{0}^{h}S_{N}(t_{l}-\sigma)\,d\sigma,\quad 1\leqslant l \leqslant m.\]
In order to have a solution in $V$, we assume that the initial value satisfies $u_0\in V$. More regularity, however, improves the spatial convergence result.
To elaborate this, we make the following assumption.
\begin{assumption}\label{assumption2a}
Let $u_0 \in\mathcal D(A^{\nu+\beta})\subset V$ for some $\beta\ge 0$. Let $\nu < 1/\rho$ and assume
that $g: [0, T] \rightarrow H : t \mapsto g(t) = f(t, u(t))$ is differentiable with bounded derivative
in $H$. Moreover, let $\gamma\ge 0$ be such that $g\in L^\infty(0,T;\mathcal D(A^\gamma))$.
\end{assumption}
Now, we are in a position to state the convergence result for the exponential Euler scheme.
\begin{theorem}\label{th.2}
 Let the initial value problem~(\ref{semilinear integro-differential equation}) satisfy Assumptions~\ref{A_operator}, \ref{smooth_property},
 \ref{assumption1}, and~\ref{assumption2a}, and consider for its numerical solution the exponential Euler method~(\ref{numerical}).
 Let $\nu<\alpha<1/\rho$.
Then, there exist constants $h_0>0$ and $C>0$ such that for all step sizes $0<h\le h_0$, the global error satisfies the bound
  \begin{eqnarray*}
  \| u(t_{m})-U^{N}_{m}\|_{V}\leqslant C \Big(t_m^{-\alpha\rho}\lambda_{N+1}^{-\alpha-\beta} + \lambda_{N+1}^{\nu-\alpha-\gamma}
  +h\sup _{0\leqslant t\leqslant T}\|g'(t)\|_{H} \Big),
  \end{eqnarray*}
  uniformly in $0\leqslant mh \leqslant T.$
\end{theorem}
\begin{proof}
 We set $g(t) = f(t, u(t))$ in~(\ref{V_O_C})
\begin{eqnarray*}
u(t_{m})=S(t_{m})u_{0}+\sum_{j=0}^{m-1}\int^{h}_{0}S(t_{m-j}-\sigma)g(t_{j}+\sigma)\,d\sigma.
\end{eqnarray*}
By using Taylor series expansion, we have
\begin{align}\label{Exact 2}
u(t_{m})=S(t_{m})u_{0}&+\sum_{j=0}^{m-1}\int^{h}_{0}S(t_{m-j}-\sigma)\,d\sigma g(t_{j})\nonumber\\
&+\sum_{j=0}^{m-1}\int^{h}_{0}S(t_{m-j}-\sigma)\int^{\sigma}_{0}g^{\prime}(t_{j}+\tau)\,d\tau \,d\sigma.
\end{align}

Let $e_m =u(t_m)-U_m^{N}$ denote the difference between
the exact and the numerical solution. By subtracting the numerical method~(\ref{numerical}) from ~(\ref{Exact 2}), and recalling~(\ref{Resolvent_operators}),~(\ref{coefficient}), we have
\begin{align*}
e_{m}&=S(t_{m})(I-\mathcal{P}_{N})u_{0} +h\sum^{m-1}_{j=0}b_{1}(t_{m-j}) (I-\mathcal{P}_{N})g(t_{j})\\
 &\quad +h\sum^{m-1}_{j=0}b_{1}^{N}(t_{m-j}) \mathcal{P}_{N} \Big(g(t_{j})-f(t_{j},U_{j}^{N})\Big)+\delta_{m},
\end{align*}
where
$$
b_{1}(t_{l})=\varphi_{1,h}(t_l)=\frac{1}{h}\int_{0}^{h}S(t_{l}-\sigma)\,d\sigma,\quad 1\leqslant l \leqslant m,
$$
and
\[ \delta_{m}=\sum_{j=0}^{m-1}\int^{h}_{0}S(t_{m-j}-\sigma)\int^{\sigma}_{0}g'(t_{j}+\tau)\,d\tau \,d\sigma.\]

By taking norms, this implies
\begin{align}\label{error1}
\|e_{m}\|_{V} &\leqslant \|S(t_{m})(I-\mathcal{P}_{N})u_{0}\|_{V} + \Big\| h\sum^{m-1}_{j=0}b_{1}(t_{m-j})
                (I-\mathcal{P}_{N})g(t_{j})\Big\|_{V}\nonumber\\
&\quad +\Big\| h\sum^{m-1}_{j=0}b_{1}^{N}(t_{m-j}) \mathcal{P}_{N} \Big(g(t_{j})-f(t_{j},U_{j}^{N})\Big)\Big\|_{V}
+ \|\delta_{m}\|_{V}=\sum_{i=1}^{4}I_{i}.
\end{align}
We note that $I_{1}$ and $I_{2}$ correspond to the spatial discretization error, while $I_{3}$  and $I_{4}$
correspond to the temporal error.

 \textit{(i) Spatial error:} The estimate of $I_{1}$ is a consequence of (\ref{property}) and
(\ref{Bound}), as
\begin{align}\label{I11}
 I_{1}= \| S(t_{m})(I-\mathcal{P}_{N})u_{0}\|_{V} &\leqslant \|A^{\alpha}S(t_{m})\|_{\mathcal{B}}\| A^{-\alpha-\beta}(I-\mathcal{P}_{N})A^{\nu + \beta}u_{0}\| \nonumber \\
                                                                          &\leqslant Ct_{m}^{-\alpha \rho}\lambda_{N+1}^{-\alpha-\beta}\Vert A^\beta u_{0}\Vert_{V}\\
                                                                          &\leqslant Ct_{m}^{-\alpha \rho}\lambda_{N+1}^{-\alpha-\beta} \nonumber.
 \end{align}
 Also for $I_{2},$ by using (\ref{property1}) and (\ref{Bound}), we have
 \begin{align}\label{I12}
 I_{2}&=\Big\| h\sum^{m-1}_{j=0}b_{1}(t_{m-j}) (I-\mathcal{P}_{N})g(t_{j})\Big\|_{V} \nonumber\\
         &\leqslant h\sum^{m-1}_{j=0}\| A^{\alpha}b_{1}(t_{m-j})\|_{\mathcal{B}}
           \| A^{\nu-\alpha-\gamma}(I-\mathcal{P}_{N})\|_{\mathcal{B}}\|A^\gamma g(t_{j})\|\\
         &\leqslant Ch\sum^{m-1}_{j=0}t_{m-j}^{-\alpha\rho}\lambda_{N+1}^{\nu-\alpha-\gamma}\|A^\gamma g(t_{j})\|\nonumber\\
         &\leqslant C\lambda_{N+1}^{\nu-\alpha-\gamma} \nonumber.
 \end{align}
  \textit{(ii) Temporal error:} Here we estimate $I_{3}$ with the help of Assumption~\ref{assumption1}, i.e.,
  \begin{align}\label{I13}
 I_{3}&=\Big\| h\sum^{m-1}_{j=0}b_{1}^{N}(t_{m-j}) \mathcal{P}_{N} \Big(g(t_{j})-f(t_{j},U_{j}^{N})\Big)\Big\|_{V} \nonumber\\
         &\leqslant h\sum^{m-1}_{j=0}\Big\| A^{\nu}b_{1}^{N}(t_{m-j}) \mathcal{P}_{N}\Big(g(t_{j})-f(t_{j},U_{j}^{N})\Big)\Big\| \\
         &\leqslant Ch\sum^{m-1}_{j=0}t_{m-j}^{-\nu\rho}\|u(t_{j})-U_{j}^{N}\|_{V}\leqslant Ch\sum^{m-1}_{j=0}t_{m-j}^{-\nu\rho}\|e_{j}\|_{V}\nonumber.
 \end{align}
Now we estimate $I_{4}$. By using~(\ref{property}), we obtain
\begin{align}\label{I14}
  \| \delta_{m}\|_{V} &\leqslant \sum_{j=0}^{m-1}\Big\| \int^{h}_{0}A^{\nu}S(t_{m-j}-\sigma)\int^{\sigma}_{0} g^{\prime}(t_{j}+\tau) \,d\tau \,d\sigma \Big\|_{H} \nonumber\\
                          & \leqslant C\sup_{0\leqslant t\leqslant T}\| g^{\prime}(t)\|_{H} \sum_{j=0}^{m-1}\int^{h}_{0}(t_{m-j}-\sigma) ^{-\nu\rho}\,\sigma \,d\sigma \nonumber\\
                          &\leqslant Ch\sup_{0\leqslant t\leqslant T}\|g^{\prime}(t)\|_{H}.
\end{align}
 Finally, inserting~(\ref{I11}),~(\ref{I12}), ~(\ref{I13}) and~(\ref{I14}) into~(\ref{error1}), we have
 \begin{align*}
 \| u(t_{m})-U^{N}_{m}\|_{V}&\leqslant Ct_m^{-\alpha\rho}\lambda_{N+1}^{-\alpha-\beta}+ C\lambda_{N+1}^{\nu-\alpha-\gamma} + Ch\sum^{m-1}_{j=0}t_{m-j}^{-\nu\rho}\|e_{j}\|_{V}+ Ch\sup_{0\leqslant t\leqslant T}\| g^{\prime}(t)\|_{H},
 \end{align*}
 which, by the discrete Gronwall lemma \cite[Lemma 2.15]{HO2}, gives
  \begin{align*}
 \| u(t_{m})-U^{N}_{m}\|_{V} \leqslant C\Big(t_m^{-\alpha\rho}\lambda_{N+1}^{-\alpha-\beta} + \lambda_{N+1}^{\nu-\alpha-\gamma}
  +h\sup _{0\leqslant t\leqslant T}\|g'(t)\|_{H} \Big).
 \end{align*}
 This is the desired result.
\end{proof}
\subsection{Convergence results for second-order methods.}
For the numerical solution of~(\ref{semilinear integro-differential equation}), we consider now second-order exponential Runge--Kutta methods, which requires two stages, i.e.~$s=2$ in~(\ref{1.14}):
\begin{subequations}\label{num_second}
\begin{align}
& U_{m}^{N}=S_{N}(t_{m})\mathcal{P}_{N}u_{0}+h\sum^{m-1}_{j=0}\sum^{2}_{i=1}b_{i}^{N}(t_{m-j}) \mathcal{P}_{N} f(t_{j}+c_{i}h,U_{j,i}^{N}),\label{num_second a}\\
& U_{m-1,1}^{N}=U_{m-1}^{N},\nonumber\\[1mm]
& U_{m-1,2}^{N}=S_{N}(t_{m-1}+c_{2}h)\mathcal{P}_{N}u_{0} + ha_{21}^{N}\mathcal{P}_{N} f(t_{m-1},U_{m-1}^{N})\nonumber\\
  & \qquad \qquad + h\sum^{m-2}_{l=0}\sum^{2}_{i=1}b_{i}^{2,N}(t_{m-l-1})\mathcal{P}_{N}f(t_{l}+c_{i}h,U_{l,i}^{N}).\label{num_second b}
\end{align}
\end{subequations}
\black Recall that we have chosen $c_1=0.$ \black

In the same way as for the exponential Euler method, we start the analysis by inserting the exact
 solution into the numerical scheme. This yields
\begin{subequations}\label{5.2}
\begin{align}
& u(t_{m})=S_{N}(t_{m})\mathcal{P}_{N}u_{0}+h\sum^{m-1}_{j=0}\sum^{2}_{i=1}b_{i}^{N}(t_{m-j}) \mathcal{P}_{N} g(t_{j}+c_{i}h)+\delta_{m},\label{5.2a}\\
& u(t_{m-1}+c_{2}h)=S_{N}(t_{m-1}+c_{2}h)\mathcal{P}_{N}u_{0} + ha_{21}^{N}\mathcal{P}_{N} g(t_{m-1})\nonumber\\
 & \qquad \qquad  +h\sum^{m-2}_{l=0}\sum^{2}_{i=1}b_{i}^{2,N}(t_{m-l-1})\mathcal{P}_{N}g(t_{l}+c_{i}h)+\Delta_{m-1,2},\label{5.2b}
\end{align}
\end{subequations}
with defects $\delta_{m}$ and $\Delta_{m-1,2}$.

Now we derive bounds for the defects $\delta_{m}$ and $\Delta_{m-1,2}$. To carry out this, we need a strengthened version of Assumption~\ref{assumption2a}.

 \begin{assumption}\label{assumption2b}
Let $u_0 \in\mathcal D(A^{\nu+\beta})\subset V$ for some $\beta\ge 0$. Let $\nu < 1/\rho$ and assume
that $g: [0, T] \rightarrow H : t \mapsto g(t) = f(t, u(t))$ is twice differentiable with bounded derivatives
in $H$. Moreover, let $\gamma\ge 0$ be such that $g\in L^\infty(0,T;\mathcal D(A^\gamma))$ and let $0\le \eta\le \nu$ be
such that $g'\in L^\infty(0,T;\mathcal D(A^\eta))$.
\end{assumption}
 By using Taylor series expansion, recalling~(\ref{Resolvent_operators}), and subtracting ~(\ref{5.2a}) from~(\ref{V_O_C}), we obtain
\begin{align*}
\delta_{m}&=S(t_{m})(I-\mathcal{P}_{N})u_{0}+\sum^{m-1}_{j=0}\int_{0}^{h}S(t_{m-j}-\sigma)(I-\mathcal{P}_{N})g(t_{j})\,d\sigma\\
&\quad +\sum^{m-1}_{j=0}\Big(\int_{0}^{h}S_{N}(t_{m-j}-\sigma)\,d\sigma-h\sum^{2}_{i=1}b_{i}^{N}(t_{m-j})\Big)\mathcal{P}_{N}g(t_{j})\\
&\quad +\sum^{m-1}_{j=0}\int_{0}^{h}S(t_{m-j}-\sigma)\sigma(I-\mathcal{P}_{N})g'(t_{j})\,d\sigma\\
&\quad +\sum^{m-1}_{j=0}\Big(\int_{0}^{h}S_{N}(t_{m-j}-\sigma)\sigma d\sigma - h^{2}b_{2}^{N}(t_{m-j})c_{2}\Big)\mathcal{P}_{N} g'(t_{j})\\
&\quad + \sum_{j=0}^{m-1}\int^{h}_{0}S(t_{m-j}-\sigma)\int^{\sigma}_{0} (\sigma - \tau)g^{\prime\prime}(t_{j}+\tau) \,d\tau \,d\sigma\\
&\quad - h\sum^{m-1}_{j=0}b_{2}^{N}(t_{m-j})\int^{c_{2}h}_{0}(c_{2}h- \tau)\mathcal{P}_{N}g^{\prime\prime}(t_{j}+\tau)\,d\tau.
\end{align*}
In order to get small defects, we choose the coefficients $b_{1}^{N}$ and $b_{2}^{N}$ such that
\begin{equation}\label{oc-part1}
\begin{aligned}
  b_{1}^{N}(t_{n}) +  b_{2}^{N}(t_{n})&=\varphi_{1,h}(t_{n})\mathcal{P}_{N},\\
  b_{2}^{N}(t_{n})c_2 &= \varphi_{2,h}(t_{n})\mathcal{P}_{N}
\end{aligned}
\end{equation}
are satisfied. These conditions are the first part of the sought-after order conditions.

In the same way, we study the stages. First, we represent the  exact solution by the variation-of-constants formula
 \begin{align}\label{Exact_tj}
u(t_{m-1}+c_{2}h)&=S(t_{m-1}+c_{2}h)u_{0}+\int_{t_{m-1}}^{t_{m-1}+c_{2}h}S(t_{m-1}+c_{2}h-\sigma)g(\sigma)\,d\sigma\nonumber\\
                         &\quad + \sum_{l=0}^{m-2}\int_{t_{l}}^{t_{l+1}}S(t_{m-1}+c_{2}h-\sigma)g(\sigma)\,d\sigma \nonumber\\
                         &=S(t_{m-1}+c_{2}h)u_{0}+\int_{0}^{c_{2}h}S(c_{2}h-\sigma)g(t_{m-1}+\sigma)\,d\sigma \nonumber\\
                         &\quad +\sum_{l=0}^{m-2}\int_{0}^{h}S(t_{m-l-1}+c_{2}h-\sigma)g(t_{l}+\sigma)\,d\sigma.
\end{align}
By using Taylor series expansion, recalling~(\ref{Resolvent_operators}), and subtracting~(\ref{5.2b}) from~(\ref{Exact_tj}), we have
\begin{align*}
\Delta_{m-1,2}&=S(t_{m-1}+c_{2}h)(I-\mathcal{P}_{N})u_{0}+\int_{0}^{c_{2}h}S(c_{2}h-\sigma)(I-\mathcal{P}_{N})g(t_{m-1})\,d\sigma\\
&\quad +\sum^{m-2}_{l=0}\int_{0}^{h}S(t_{m-l-1}+c_{2}h-\sigma)(I-\mathcal{P}_{N})g(t_{l})\,d\sigma\\
&\quad +\sum^{m-2}_{l=0}\int_{0}^{h}S(t_{m-l-1}+c_{2}h-\sigma)\sigma(I-\mathcal{P}_{N})g'(t_{l})\,d\sigma\\
&\quad +\Big(\int_{0}^{c_{2}h}S_{N}(c_{2}h-\sigma)d\sigma -ha_{21}^{N}\Big) \mathcal{P}_{N}g(t_{m-1})\\
&\quad +\int^{c_{2}h}_{0}S(c_{2}h-\sigma)\int^{\sigma}_{0} g'(t_{m-1}+\tau)\,d\tau \,d\sigma\\
&\quad +\sum^{m-2}_{j=0}\Big(\int_{0}^{h}S_{N}(t_{m-l-1}+c_{2}h-\sigma)d\sigma -h\sum^{2}_{i=1}b_{i}^{2,N}(t_{m-l-1})\Big)
\mathcal{P}_{N}g(t_{l})\\
&\quad +\sum^{m-2}_{l=0}\Big(\int_{0}^{h}S_{N}(t_{m-l-1}+c_{2}h-\sigma)\sigma \,d\sigma- h^{2}b_{2}^{2,N}(t_{m-l-1})c_{2}\Big) \mathcal{P}_{N}g'(t_{l})\\
&\quad + \sum_{l=0}^{m-2}\int^{h}_{0}S(t_{m-l-1}+c_{2}h-\sigma)\int^{\sigma}_{0} (\sigma - \tau)g^{\prime\prime}(t_{l}+\tau) \,d\tau \,d\sigma\\
&\quad - h\sum^{m-2}_{l=0}b_{2}^{2,N}(t_{m-l-1})\int^{c_{2}h}_{0}(c_{2}h- \tau)\mathcal{P}_{N}g^{\prime\prime}(t_{l}+\tau)\,d\tau.
\end{align*}

Again, the coefficients are chosen to minimize the defects. This results in
\begin{equation}\label{oc-part2}
\begin{aligned}
  a_{21}^N &= c_2 \varphi_{1,c_2h}(c_2h)\mathcal{P}_{N} = \frac1h\int_0^{c_2h} S_N(c_2 h-\sigma)\,d\sigma,\\
  b_{1}^{2,N}(t_{n}) +  b_{2}^{2,N}(t_{n})&=\varphi_{1,h}(t_n+c_2h)\mathcal{P}_{N},\\
  b_{2}^{2,N}(t_{n})c_2 &= \varphi_{2,h}(t_{n}+c_2h)\mathcal{P}_{N},
\end{aligned}
\end{equation}
which is the second set of order conditions. The final set of order conditions of order two is given in Table~\ref{table2}.

 \begin{table}[ht]
 \caption{Order conditions \eqref{oc-part1} and \eqref{oc-part2} for a two-stage explicit exponential Runge--Kutta
 methods applied to~\eqref{semilinear integro-differential equation}. The functions $\varphi_{k,h}$ are defined in~(\ref{phi}).\label{table2}}
 \renewcommand{\arraystretch}{1.8}
 \begin{tabular}{|c|c|l|}
 \hline
$\text{Number}$ &$\text{Order}$ & $\qquad \qquad \qquad\text{Order condition}$ \\
 \hline
 $1$ &$1$ & $ b_1^N(t_n) + b_2^N(t_n)= \varphi_{1,h}(t_{n})\mathcal{P}_{N},\quad t_n\in[0,T]$ \\
\hline
 $2$ &$2$ & $ b_{2}^N(t_{n})c_{2}=\varphi_{2,h}(t_{n})\mathcal{P}_{N},\quad t_n\in[0,T]$ \\
 $3$ &$2$ & $ a_{21}^{N}= c_2\varphi_{1,c_2 h}(c_2 h)\mathcal{P}_{N}$\\
 $4$ &$2$ & $b_{1}^{2,N}(t_{l}) + b_{2}^{2,N}(t_{l})=\varphi_{1,h}(t_{l}+c_{2}h)\mathcal{P}_{N},\quad t_l+c_2h\in[0,T]$ \\
 $5$ &$2$ & $ b^{2,N}_{2}(t_{l})c_{2}=\varphi_{2,h}(t_{l}+c_{2}h)\mathcal{P}_{N},\quad t_l+c_2h\in[0,T]$ \\
\hline
 \end{tabular}
\end{table}
Using the order conditions of Table 2, we can derive bounds for the defects $\delta_{m}$ and $\Delta_{m-1,2}$. By taking the norm of $\delta_{m}$, we have
\begin{align*}
\| \delta_{m}\|_{V}&\leqslant \| S(t_{m})(I-\mathcal{P}_{N})u_{0}\|_{V} +\| \sum^{m-1}_{j=0}\int_{0}^{h}S(t_{m-j}-\sigma)(I-\mathcal{P}_{N})g(t_{j})\,d\sigma\|_{V}\\
&\quad +\|\sum^{m-1}_{j=0}\int_{0}^{h}S(t_{m-j}-\sigma)\sigma(I-\mathcal{P}_{N})g'(t_{j})\,d\sigma\|_{V}\\
&\quad +\| \sum_{j=0}^{m-1}\int^{h}_{0}S(t_{m-j}-\sigma)\int^{\sigma}_{0} (\sigma - \tau)g^{\prime\prime}(t_{j}+\tau) \,d\tau \,d\sigma\|_{V}\\
&\quad +\| h\sum^{m-1}_{j=0}b_{2}^{N}(t_{m-j})\int^{c_{2}h}_{0}(c_{2}h- \tau)\mathcal{P}_{N}g^{\prime\prime}(t_{j}+\tau)\,d\tau\|_{V}=\sum_{i=1}^{5}\delta_{mi}.
\end{align*}
We note that $\delta_{m1}$, $\delta_{m2}$, and $\delta_{m3}$ correspond to the spatial discretization error. Under Assumption~\ref{assumption2b} these terms are estimated in the same way as the corresponding terms for the exponential Euler scheme. Therefore, we have
\begin{align}\label{delta}
  \| \delta_{m}\|_{V} &\leqslant  Ct_m^{-\alpha\rho}\lambda_{N+1}^{-\alpha-\beta}+ C\lambda_{N+1}^{\nu-\alpha-\gamma} \nonumber \\
  &\quad +\sum_{j=0}^{m-1}\Big\| \int^{h}_{0}A^{\nu}S(t_{m-j}-\sigma)\int^{\sigma}_{0} (\sigma - \tau)g^{\prime\prime}(t_{j}+\tau) \,d\tau \,d\sigma \Big\|_{H} \nonumber\\
 &\quad +h\sum^{m-1}_{j=0}\Big\|A^{\nu}b_{2}^{N}(t_{m-j})\int^{c_{2}h}_{0}(c_{2}h- \tau)\mathcal{P}_{N}g^{\prime\prime}(t_{j}+\tau)\,d\tau \Big\|_{H} \nonumber\\
  &\leqslant  Ct_m^{-\alpha\rho}\lambda_{N+1}^{-\alpha-\beta}+ C\lambda_{N+1}^{\nu-\alpha-\gamma} \nonumber\\
  &\quad +C\sup_{0\leqslant t\leqslant T}\| g^{\prime\prime}(t)\|_{H} \sum_{j=0}^{m-1}\int^{h}_{0}(t_{m-j}-\sigma)^{-\nu\rho} \int^{\sigma}_{0}(\sigma - \tau)\,d\tau \,d\sigma \nonumber\\
  & \quad +C\sup_{0\leqslant t\leqslant T}\| g^{\prime\prime}(t)\|_{H}\cdot h\sum_{j=0}^{m-1}t_{m-j}^ {-\nu\rho} \int^{c_{2}h}_{0}(c_{2}h - \tau)\,d\tau \nonumber\\
  &\leqslant Ct_m^{-\alpha\rho}\lambda_{N+1}^{-\alpha-\beta}+ C\lambda_{N+1}^{\nu-\alpha-\gamma}+Ch^{2}\sup_{0\leqslant t\leqslant T}\| g^{\prime\prime}(t)\|_{H}.
\end{align}
Also, by taking the norm of $\Delta_{m-1,2}$, we have
\begin{align*}
\|\Delta_{m-1,2}\|_{V} &\leqslant \|S(t_{m-1}+c_{2}h)(I-\mathcal{P}_{N})u_{0}\|_{V}+\Big\|\int_{0}^{c_{2}h}S(c_{2}h-\sigma)(I-\mathcal{P}_{N})g(t_{m-1})\,d\sigma\Big\|_{V}\\
&\quad +\Big\|\sum^{m-2}_{l=0}\int_{0}^{h}S(t_{m-l-1}+c_{2}h-\sigma)(I-\mathcal{P}_{N})g(t_{l})\,d\sigma\Big\|_{V}\\
&\quad +\Big\|\sum^{m-2}_{l=0}\int_{0}^{h}S(t_{m-l-1}+c_{2}h-\sigma)\sigma(I-\mathcal{P}_{N})g'(t_{l})\,d\sigma\Big\|_{V}\\
&\quad  +\Big\|\int^{c_{2}h}_{0}S(c_{2}h-\sigma) \int^{\sigma}_{0}g'(t_{j}+\tau)\,d\tau \,d\sigma\Big\|_{V}\\
&\quad +\Big\| \sum_{l=0}^{m-2}\int^{h}_{0}S(t_{m-l-1}+c_{2}h-\sigma)\int^{\sigma}_{0} (\sigma - \tau)g^{\prime\prime}(t_{l}+\tau) \,d\tau \,d\sigma\Big\|_{V}\\
&\quad +\Big\| h\sum^{m-2}_{l=0}b_{2}^{2,N}(t_{m-l-1})\int^{c_{2}h}_{0}(c_{2}h- \tau)\mathcal{P}_{N}g^{\prime\prime}(t_{l}+\tau)\,d\tau\Big\|_{V}=\sum_{i=1}^{7}\|\Delta_{m-1,2}^{i}\|.
\end{align*}
The terms $\Delta_{m-1,2}^{1}$ to $\Delta_{m-1,2}^{4}$ correspond to the spatial discretization error, so we get
  \begin{align*}
\|\Delta_{m-1,2}\|_{V} &\leqslant  Ct_m^{-\alpha\rho}\lambda_{N+1}^{-\alpha-\beta}+ C\lambda_{N+1}^{\nu-\alpha-\gamma} \nonumber \\
           &\quad +C\sup _{0\leqslant t\leqslant T} \|A^\eta g'(t)\|_{H} \int_{0}^{c_{2}h}(c_{2}h-\sigma)^{-(\nu-\eta)\rho} \sigma \,d\sigma\\
           &\quad +C\sup _{0\leqslant t\leqslant T}\| g^{\prime\prime}(t)\|_{H} \sum_{l=0}^{m-2}\int_{0}^{h}(t_{m-l-1}+c_{2}h-\sigma )^{-\nu\rho}\sigma^{2}\,d\sigma\\
           &\quad +Ch^{2}\sup _{0\leqslant t\leqslant T}\| g^{\prime\prime}(t)\|_{H} \cdot h\sum_{l=0}^{m-2}t_{m-l-1} ^{-\nu\rho}
 \end{align*}
and finally
 \begin{align*}
\|\Delta_{m-1,2}\|_{V} &\leqslant  Ct_m^{-\alpha\rho}\lambda_{N+1}^{-\alpha-\beta}+ C\lambda_{N+1}^{\nu-\alpha-\gamma} \\
           &\quad + Ch^{2-(\nu-\eta)\rho}\sup _{0\leqslant t\leqslant T}\| A^\eta g'(t)\|_{H} +Ch^{2}\sup _{0\leqslant t\leqslant T}\| g^{\prime\prime}(t)\|_{H}.
  \end{align*}
  Now we are ready to state our convergence result.
\begin{theorem}\label{th.3}
Let the initial value problem~(\ref{semilinear integro-differential equation}) satisfy Assumptions~\ref{A_operator}, \ref{smooth_property},
 \ref{assumption1}, and~\ref{assumption2b}, and consider for its numerical solution the exponential Runge--Kutta method (\ref{num_second}) that satisfies the order conditions of Table~\ref{table2}. Let $\nu<\alpha<1/\rho$.
Then, there exist constants $h_0>0$ and $C>0$ such that for all step sizes $0<h\le h_0$, the global error satisfies the bound
  \begin{align*}
  \| u(t_{m})-U^{N}_{m}\|_{V}&\leqslant C\Big(t_m^{-\alpha\rho}\lambda_{N+1}^{-\alpha-\beta} + \lambda_{N+1}^{\nu-\alpha-\gamma}\\
  &\qquad\qquad +h^{2-(\nu-\eta)\rho}\sup _{0\leqslant t\leqslant T}\|A^\eta g'(t)\|_{H}  +h^2\sup _{0\leqslant t\leqslant T}\|g''(t)\|_{H} \Big),
  \end{align*}
  uniformly in $0\leqslant mh \leqslant T.$
\end{theorem}
In particular, if $g'$ is uniformly bounded in $V$, we can choose $\eta=\nu$ and the scheme turns out to be second-order convergent in time.
\begin{proof}
Let $e_m =u(t_m)-U_m^{N}$ and $E_{j,2}=u(t_{j} + c_{2}h)- U_{j,2}^{N}$ denote the differences between
the exact solution~(\ref{5.2}) and the numerical solution~(\ref{num_second}). Then
\begin{align*}
e_{m}=h\sum^{m-1}_{j=0}\sum^{2}_{i=1}b_{i}^{N}(t_{m-j}) \mathcal{P}_{N} \Big(g(t_{j}+c_{i}h)-f(t_{j}+c_{i}h,U_{j,i}^{N})\Big) +\delta_{m}.
\end{align*}
By taking norms, we obtain
\begin{align}\label{error}
\|e_{m}\|_{V} &\leqslant \Big\| h\sum^{m-1}_{j=0}b_{1}^{N}(t_{m-j}) \mathcal{P}_{N} \Big(g(t_{j})-f(t_{j},U_{j}^{N})\Big)\Big\|_{V} \nonumber\\
 &\quad +\Big\| h\sum^{m-1}_{j=0}b_{2}^{N}(t_{m-j}) \mathcal{P}_{N} \Big(g(t_{j}+c_{2}h)-f(t_{j}+c_{2}h,U_{j,2}^{N})\Big)\Big\|_{V}+ \|\delta_{m}\|_{V}\nonumber\\
 &=\sum_{i=1}^{3}I_{i}.
\end{align}
We know that $I_{1}$ and $I_{2}$ correspond to the temporal error. The term $I_{3}$ has already been estimated.

 First we bound $I_{1}$, i.e.,
 \begin{align}\label{I3}
 I_{1}&=\Big\| h\sum^{m-1}_{j=0}b_{1}^{N}(t_{m-j}) \mathcal{P}_{N} \Big(g(t_{j})-f(t_{j},U_{j}^{N})\Big)\Big\|_{V} \nonumber\\
         &\leqslant h\sum^{m-1}_{j=0}\Big\| A^{\nu} b_{1}^{N}(t_{m-j})\mathcal{P}_{N}\Big(g(t_{j})-f(t_{j},U_{j}^{N})\Big)\Big\| \\
         &\leqslant Ch\sum^{m-1}_{j=0}t_{m-j}^{-\nu\rho}\| u(t_{j})-U_{j}^{N}\|_{V}= Ch\sum^{m-1}_{j=0}t_{m-j}^{-\nu\rho}\|e_{j}\|_{V}\nonumber.
 \end{align}
 Now  we estimate $I_{2}$,
 \begin{align}\label{I4}
 I_{2}&=\Big\| h\sum^{m-1}_{j=0}b_{2}^{N}(t_{m-j}) \mathcal{P}_{N} \Big(g(t_{j}+c_{2}h)-f(t_{j}+c_{2}h,U_{j,2}^{N})\Big)\Big\|_{V} \nonumber\\
        &\leqslant h\sum^{m-1}_{j=0}\Big\| A^{\nu}b_{2}^{N}(t_{m-j})\mathcal{P}_{N}\Big(g(t_{j}+c_{2}h)-f(t_{j}+c_{2}k,U_{j,2}^{N})\Big)\Big\|\nonumber\\
         &\leqslant Ch\sum^{m-1}_{j=0}t_{m-j}^{-\nu\rho}\| g(t_{j}+c_{2}h)-f(t_{j}+c_{2}h,U_{j,2}^{N})\|\nonumber\\
         &\leqslant Ch\sum^{m-1}_{j=0} t_{m-j}^{-\nu\rho}\|E_{j,2}\|_{V}.
 \end{align}
  For $E_{m-1,2}=u(t_{m-1} + c_{2}h)- U_{m-1,2}^{N}$, we have
 \begin{align*}
 E_{m-1,2}&= ha_{21}^{N}\mathcal{P}_{N}\Big(g(t_{m-1})-f(t_{m-1},U_{m-1}^{N})\Big)
+ h\sum^{m-2}_{l=0}b_{1}^{2,N}(t_{m-1-l})\mathcal{P}_{N}\Big(g(t_{l})-f(t_{l},U_{l}^{N})\Big)\\
  &\quad + h\sum^{m-2}_{l=0}b_{2}^{2,N}(t_{m-1-l})\mathcal{P}_{N}\Big(g(t_{l}+c_{2}h)-f(t_{l}+c_{2}h,U_{l,2}^{N})\Big)+\Delta_{m-1,2}.
 \end{align*}
 By taking norm, we get
 \begin{align*}
 \| E_{m-1,2}\|_{V} &\leqslant \Big \| ha_{21}^{N}\mathcal{P}_{N}\Big(g(t_{m-1})-f(t_{m-1},U_{m-1}^{N})\Big)\Big\|_{V} \nonumber\\
 &\quad +\Big\| h\sum^{m-2}_{l=0}b_{1}^{2,N}(t_{m-1-l})\mathcal{P}_{N}\Big(g(t_{l})-f(t_{l},U_{l}^{N})\Big)\Big\|_{V} \\
 &\quad + \Big\| h\sum^{m-2}_{l=0}b_{2}^{2,N}(t_{m-1-l})\mathcal{P}_{N}\Big(g(t_{l}+c_{2}h)-f(t_{l}+c_{2}h,U_{l,2}^{N})\Big)\Big\|_{V} +\|\Delta_{m-1,2}\|_{V}\\
 & =\sum^{4}_{j=1}I_{2,j}.
 \end{align*}
 For $I_{2,1}$, $I_{2,2}$ and $I_{2,3}$, we have
 \begin{align}\label{I456}
 I_{2,1}+I_{2,2}+I_{2,3}\leqslant Ch^{1-\nu\rho} \|e_{m-1}\|_{V} + Ch\sum^{m-2}_{j=0}t_{m-1-j}^{-\nu\rho}\|e_{j}\|_{V} +
 Ch\sum^{m-2}_{j=0}t_{m-1-j}^{-\nu\rho}\|E_{j,2}\|_{V}.
 \end{align}
 Taking all together, we obtain
 \begin{align*}
 \|e_m\|_V \le Ch\sum^{m-1}_{j=0}t_{m-j}^{-\nu\rho}\bigl(\|e_{j}\|_{V}+ \|E_{j,2}\|\bigr) + \|\delta_m\|
 \end{align*}
 and
 \begin{align*}
\|E_{m-1,2}\|_V \le Ch^{1-\nu\rho} \|e_{m-1}\|_{V} + Ch\sum^{m-2}_{j=0}t_{m-1-j}^{-\nu\rho}\bigl(\|e_{j}\|_{V}+ \|E_{j,2}\|\bigr) + \|\Delta_{m-1,2}\|.
 \end{align*}
Applying a discrete Gronwall lemma~\cite[Lemma 2.15]{HO2} finally gives the desired result.
\end{proof}
\section{Numerical implementation.}\label{Numerical}
In this section we first derive an explicit representation of the resolvent family for two different kernels for the problem
\begin{eqnarray*}
 u'(t) +\int_{0}^{t}b(t-s)Au(s)\,ds=f(t,u(t)), \quad t\in (0,T], \quad  u(0)=u_{0}.
\end{eqnarray*}
 Then, we illustrate by numerical experiments the temporal order of convergence, to confirm the rates proposed in
 Theorems \ref{th.2} and \ref{th.3}.

Let $\lbrace(\lambda_{k},\psi_{k})\rbrace_{k=1}^{\infty}$ be the eigenpairs of $A$, i.e.,
\begin{equation*}
A\psi_{k}=\lambda_{k}\psi_{k},\quad k \in \mathbb{N}.
\end{equation*}
 Then, the resolvent family is given by
\begin{eqnarray*}
S(t)v=\sum_{k=1}^{\infty}s_{k}(t)(v,\psi_{k})\psi_{k}.
\end{eqnarray*}

To explain the implementation of the fully discrete methods, (\ref{numerical}) and~(\ref{num_second}),
we note that
\begin{eqnarray*}
S_{N}(t_{m})v=\sum_{k=1}^{N}s_{k}(t)(v,\psi_{k})\psi_{k}.
\end{eqnarray*}
The functions ${s}_{k}$ are the solutions of the scalar problems
\begin{eqnarray*}\label{equation-2}
s'_{k}(t)+\lambda_{k}\int_{0}^{t}b(t-s)s_{k}(s)\,ds=0,\quad t>0,\quad s_{k}(0)=1.
\end{eqnarray*}
In the following examples, we consider two different kernels: a Riesz kernel and an exponential kernel. We
choose
\begin{eqnarray}
\black A = -\frac{\partial^2}{\partial x^2},\black \qquad \Omega = (0,1)\subset \R.
\end{eqnarray}
For this choice, we have $\psi_{k}(x)=\sqrt{2}\sin k\pi x$, $\lambda_{k}=k^{2}\pi^{2}$ for $x\in\Omega$ and every $k \in \N$.
\begin{example}\label{Riesz}
Let $b$ be the Riesz kernel, given by $b(t)=\dfrac{t^{\beta-1}}{\Gamma(\beta)}$ for some $0<\beta<1$. We denote henceforth
$$\rho=\beta+1,\qquad 1<\rho<2,$$
so that $b(t)=\dfrac{t^{\rho-2}}{\Gamma(\rho -1)}.$
 By taking the Laplace transform of \eqref{equation-2}, we have
\begin{align*}
s_{k}(t)=E_{\rho}(-\lambda_{k}t^{\rho}),
\end{align*}
where  $E_{\rho}(-\lambda_{k}t^{\rho})$ is the one-parameter Mittag-Leffler function. Thus the resolvent family is given by
\begin{eqnarray*}
S_{N}(t)v=\sum_{k=1}^{N}E_{\rho}(-\lambda_{k}t^{\rho})(v,\psi_{k})\psi_{k}.
\end{eqnarray*}
We note that integrals of the Mittag-Leffler functions are easily computable, e.g.,
by means of a simple quadrature. The integral can be even computed exactly as
\begin{align*}
 \int_{t_{j}}^{t_{j+1}} E_{\rho}(-\lambda_{k}(t_{m}-\sigma)^{\rho})\,d\sigma &=\int_{t_{m-j-1}}^{t_{m-j}} E_{\rho}(-\lambda_{k}\sigma^{\rho})\,d\sigma\\
 &=E_{\rho,2}(-\lambda_{k}t_{m-j}^{\rho})-E_{\rho,2}(-\lambda_{k}t_{m-j-1}^{\rho}),
 \end{align*}
see \cite[Equation (1.100)]{Podlubny}.
 For evaluating the Mittag-Leffler function we use the model function from \cite{PK}.
\end{example}
 \begin{example}\label{Exponential}
 Let the kernel $b$ be an exponential function, $b(t)=e^{-at}$ with $0<a\leqslant 2$. By taking the Laplace transform of \eqref{equation}, a simple calculation shows that
\begin{align*}
s_{k}(t)=e^{\frac{-a}{2}t}\Big\{\cos\sqrt{\frac{4\lambda_{k}-a^{2}}{4}}t+\frac{a}{\sqrt{4\lambda_{k}-a^{2}}}\sin\sqrt{\frac{4\lambda_{k}-a^{2}}{4}}t\Big\}.
\end{align*}
 In our numerical experiments, we will take $a=2$.
\end{example}

\subsection{Numerical experiments.}

  We carry out experiments for the exponential Euler integrator (\ref{numerical}) and an exponential Runge--Kutta method of order two.
 Using the order conditions of Table~\ref{table2}, the coefficients of the second-order method are uniquely defined in terms of the node $c_2$:
 \begin{equation}\label{co2nd}
 \begin{aligned}
  b_{1}^{N}(t_{n})&=\varphi_{1,h}(t_{n})\mathcal{P}_{N}-\tfrac{1}{c_{2}}\varphi_{2,h}(t_{n})\mathcal{P}_{N},\\
  b_{2}^{N}(t_{n})&=\tfrac{1}{c_{2}}\varphi_{2,h}(t_{n})\mathcal{P}_{N},\\
  a_{21}^{N}&=c_2\varphi_{1,c_{2}h}(c_{2}h)\mathcal{P}_{N},\\
  b_{1}^{2,N}(t_{l})&=\varphi_{1,h}(t_{l}+c_{2}h)\mathcal{P}_{N}-\tfrac{1}{c_{2}}\varphi_{2,h}(t_{l}+c_{2}h)\mathcal{P}_{N},\\
  b_{2}^{2,N}(t_{l})&=\tfrac{1}{c_{2}}\varphi_{2,h}(t_{l}+c_{2}h)\mathcal{P}_{N}.
  \end{aligned}
\end{equation}
In our experiments, we have chosen $c_2=\frac12$.

As example, we consider \black $A = -\frac{\partial^2}{\partial x^2}$ \black on $\Omega = (0,1)$ with initial data $u_{0}=\sin(\pi x)/\sqrt{2}$ and nonlinearity $f(u(x,t))=\sin(u(x,t))$ for $x\in [0,1]$ and $t \in [0,1]$, subject to homogeneous Dirichlet boundary conditions. \black We determine a reference solution by using a very small time step (half of the smallest time step that we consider for the numerical solutions). The error is then calculated as the $L_{2}$-norm of the difference between the solution at larger time steps and the reference solution, obtained with the small time step.\black

\begin{figure}[h]
  \includegraphics[scale=0.5]{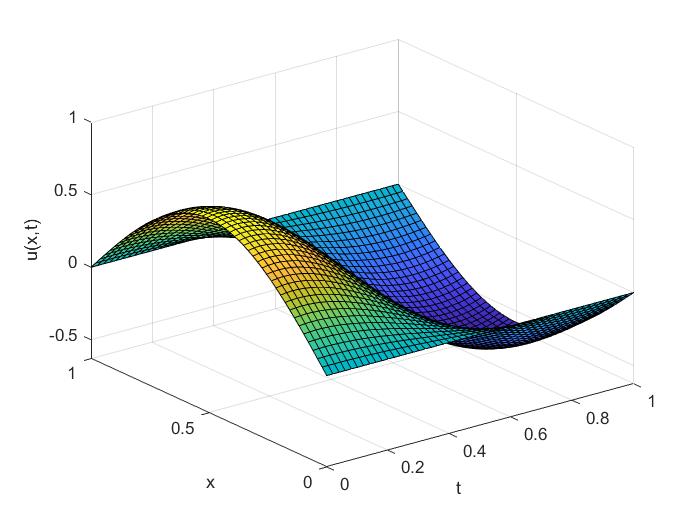}
 \caption{\small Behavior of the solution for the Riesz kernel with $\rho=1.75$.\label{fig1}}
 \end{figure}

\begin{figure}[h]
 \includegraphics[scale=0.5]{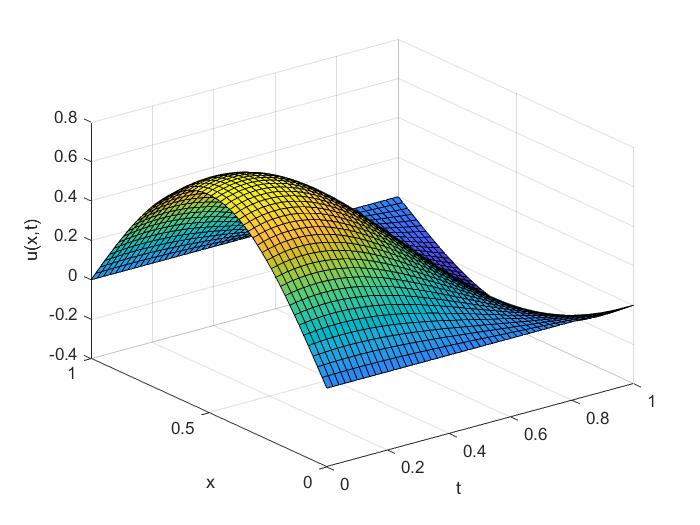}
 \caption{\small  Behavior of the solution for the exponential kernel with $a=2$.\label{fig2}}
\end{figure}

\begin{figure}[!ht]
\centering
\begin{subfigure}{0.45\textwidth}
\centering
\includegraphics[width=1.1\linewidth]{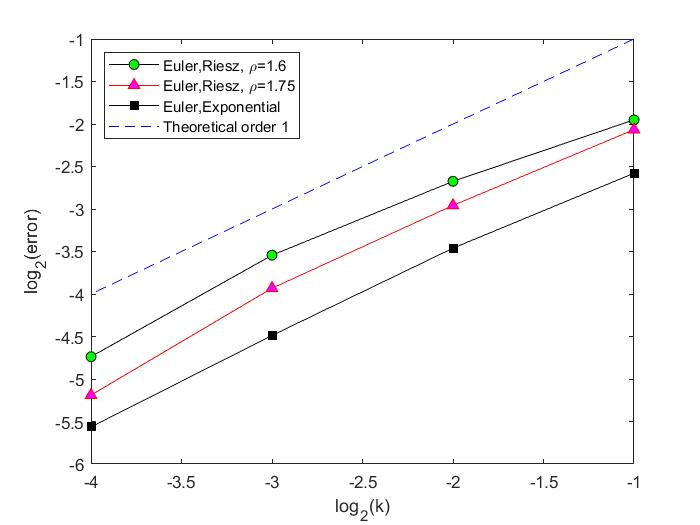}
\label{fig:sub1}
\end{subfigure}
\begin{subfigure}{0.45\textwidth}
\centering
\includegraphics[width=1.1\linewidth]{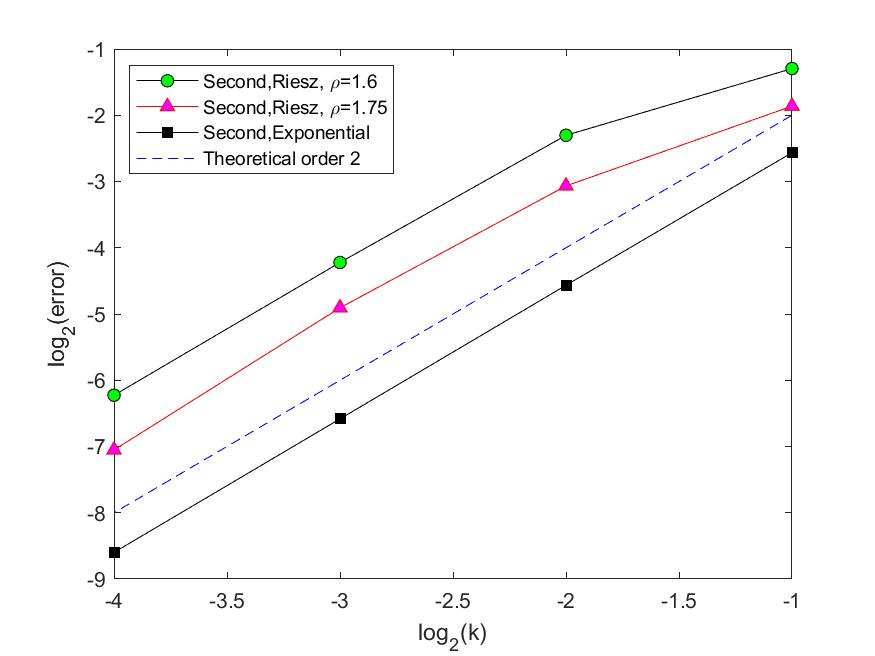}
\label{fig:sub2}
\end{subfigure}
\caption{\small  In the left panel, the temporal rate of convergence of the exponential Euler integrator applied to examples~\ref{Riesz} and~\ref{Exponential} is shown. In the right panel, the temporal rate of convergence of the second-order method is shown. \label{fig3}}
\end{figure}

We discretize this example in space by the spectral Galerkin method with 2500 points. Due to our theory, we expect to see order one for~(\ref{numerical}) with the coefficient $b_{1}^{N} =\varphi_{1,h}\mathcal{P}_{N},$
and order two for ~(\ref{num_second}) with the coefficients~\eqref{co2nd} and $c_{2} = \frac{1}{2}$. We consider different values $\rho$ for the Riesz kernel and $a=2$ for  the exponential kernel.  Figures~\ref{fig1} and~\ref{fig2} display the behaviour of the solutions. The stated orders of convergence are confirmed in Figure~\ref{fig3}. \\

\textbf{Acknowledgments.}
The main part of this work was carried out during a visit of N.~Vaisi~at the University of Innsbruck, Austria. She wishes to acknowledge their hospitality. Her research stay was financially supported by the University of Kurdistan, Iran.
\vskip 25pt
\bibliographystyle{model1-num-names}

\end{document}